\documentclass[12pt]{article}
 \usepackage[german]{babel}
\usepackage{amssymb}
\usepackage{amsmath}
\usepackage{graphicx}
\usepackage{bbm}
\usepackage[right]{eurosym}

\def\vs{\vspace}
\def\noi{\noindent}

\def\IN{\mathbb N}

\def\IR{\mathbb R}

\def\IQ{\mathbb Q}

\def\an{\mathrm{an}}
\def\exp{\mathrm{exp}}

\def\ma{\mathcal}

\begin{document}
\begin{center}
{\bf \Large Differentiability Properties of 
	
	Log-Analytic Functions}
\end{center}

\centerline{Tobias Kaiser and Andre Opris}

\vspace{0.7cm}\noi \footnotesize {{\bf Abstract.}
	We show that the derivative of a log-analytic function is log-analytic. We prove that log-analytic functions exhibit strong quasianalytic properties.
	We establish the parametric version of Tamm's theorem for log-analytic functions.

\normalsize
\section*{Introduction}

This paper contributes to analysis in the framework of o-minimal structures. O-minimality is a concept from mathematical logic with connections and applications to geometry, analysis, number theory and other areas. Sets and functions definable in an o-minimal structure (i.e. `belonging to') exhibit tame geometric and combinatorial behaviour. We refer to the book of Van den Dries [2] or to Miller and Van den Dries [5] for the general properties of o-minimal structures; in the preliminary section we state the definition and give examples. 

\;\;\;\;\;\;\;\; We consider log-analytic functions which have been defined by Lion and Rolin in their seminal paper [8]. 
They are iterated compositions from either side of globally subanalytic functions (see [5]) and the global logarithm. Their definition is kind of hybrid as we will explain. The globally subanalytic functions are precisely the functions that are definable in the o-minimal structure $\IR_\an$ of restricted analytic functions (see [5]). Since the global logarithm is not globally subanalytic the class of log-analytic functions contains properly the class of $\IR_\an$-definable functions. Observing that the exponential function is definable from the global logarithm and is not log-analytic we see that the class of log-analytic functions is not a class of definable functions. It is properly contained in the class of $\IR_{\an,\exp}$-definable functions ($\IR_{\an,\exp}$ is the o-minimal structure of restricted analytic function and the global exponential function (see [5])). Hence log-analytic functions capture $\IR_\an$-definability but not full definability; with respect to the global logarithm we restrict to composition.

\;\;\;\;\;\;\;\; Log-analytic functions are interesting from the view point of analysis. The globally subanalytic functions are not stable under parametric integration (the antiderivative of the reciprocal being the logarithm). But parametric integrals of globally subanalytic functions are log-analytic (see Lion and Rolin [9]; see also Cluckers and Dan Miller [1]).

\vs{0.5cm}
\hrule

\vs{0.4cm}
{\footnotesize{\itshape 2010 Mathematics Subject Classification:} 03C64, 14P15, 26A09, 26E05, 26E10, 32B20}
\newline
{\footnotesize{\itshape Keywords and phrases:} log-analytic functions, preparation, differentiability, Tamm's theorem}

\newpage

Globally subanalytic functions are `close' to analytic functions (see Van den Dries and Miller [4]).  They do not exhibit properties of the function $\exp(-1/x^2)$ as flatness or infinite differentiability but not real analyticity. Since the definition of log-analytic functions avoids the exponential function they should be also `close' to analytic functions.
This seems to be obvious. But the problem is that a composition of globally subanalytic functions and the logarithm allows a representation by `nice' terms only piecewise. Moreover the `pieces' are in general not definable in $\IR_\an$ but only in $\IR_{\an,\exp}$ as the following example shows (where the global logarithm is extended by $0$ to the real line).

\vs{0.5cm}
{\bf Example}

\vs{0.1cm}
Let $F(x,y)$ be the globally subanalytic function with $F(x,y)=y-x$ if $y>x$ and zero else.
Let $f(x,y)$ be the log-analytic function $F(x,\log(y))$. Then
$f(x,y)=\log(y)-x$ if $y>\exp(x)$ and zero else.

\vs{0.5cm}
If a log-analytic function has on some `piece' a representation by nice terms then on this piece the function is real analytic. 
But it is much harder to control the function on the boundary of the piece as the following example indicates.

\vs{0.5cm}
{\bf Example}

\vs{0.1cm}
Let $f(x,y)$ be the log-analytic function with
$f(x,y)=y-x/\log(y)$ if $y>0$ and zero else. Then the following asymptotics hold.
For every $x\neq 0$ we have $f(x,-)\sim -1/\log(y)$ and $f(0,-)\sim y$ as $y\searrow 0$.

\vs{0.5cm}
Hence one needs a better representation than a piecewise description by terms.
Such a good representation is given by the preparation theorems of Lion and Rolin in [8] (see also the corrections by Paw\l ucki and Pi\c{e}kosz [11]; see also Van den Dries and Speissegger [6] for a weaker preparation in a more general context; compare also [7] for another kind of preparation for the special class of constructible functions not in terms of units but suitable for questions on integration).
In case of log-analytic functions it states that the log-analytic function $f(x,y)$ where $y$ is the last variable can be piecewise written  
as 
$f(x,y)=a(x)y_0(x,y)^{q_0}\cdots y_r(x,y)^{q_r}u(x,y)$ where $y_0(x,y)=y-\Theta_0(x),y_1(x,y)=\log(|y_0(x,y)|)-\Theta_1(x),\ldots$, the $q_j$'s are rational exponents and $u(x,y)$ is a unit of a special form.
This gives roughly that the functions $f(x,-)$ behave piecewise as iterated logarithms independently of $x$ where the order of iteration is bounded in terms of $f$.
But the problem is that the functions $a(x),\Theta_0(x),\ldots,\Theta_r(x)$ although being definable in $\IR_{\an,\exp}$ are in general not log-analytic anymore. We will present an example below.
To be able to use the preparation theorems our initial result is the key observation that on certain 'pieces' which we call simple a log-analytic function can be prepared with log-analytic data only. For this one has to redo the proof of the existing preparation result.
This allows us to prove that the class of log-analytic functions is closed under taking derivatives.

\vs{0.5cm}
{\bf Theorem A}

\vs{0.1cm}
{\it Let $U\subset \IR^n$ be open and let $f:U\to \IR$ be log-analytic. Let $i\in \{1,\ldots,m\}$ be such that $f$ is differentiable with respect to $x_i$ on $U$.  Then $\partial f/\partial x_i$ is log-analytic.}

\vs{0.5cm}
The asymptotic behaviour of log-analytic functions given by the preparation results on simple sets implies strong quasianalyticity (see Miller [10] for this result in polynomially bounded o-minimal structures). 

\vs{0.5cm}
{\bf Theorem B}

\vs{0.1cm}
{\it Let $U\subset \IR^n$ be a domain and let $f:U\to \IR$ be a log-analytic function. Then there is $N\in \IN$ with the following property. If $f$ is $C^N$ and if there is $a\in U$ such that all derivaties of $f$ up to order $N$ vanish in $a$ then $f$ vanishes identically.}

\vs{0.5cm}
With the results above, we can generalize the parametric version of Tamm's theorem for globally subanalytic functions by Van den Dries and Miller [4] (see Tamm [12] for the original version) to log-analytic functions by adapting their arguments.

\vs{0.5cm}
{\bf Theorem C} 

\vs{0.1cm}
{\it Let $f:\IR^n\times\IR^m\to \IR, (x,y)\mapsto f(x,y),$ be a log-analytic function.
	Then there is $N\in \IN$ such that the following holds for every $(x,y)\in \IR^n\times\IR^m$.
	If $f(x,-)$ is $C^N$ at $y$ then $f(x,-)$ is real analytic at $y$.} 

\vs{0.5cm}
This implies in particular that the set of all $x\in \IR^n$ such that $f(x,-)$ is real analytic is definable in $\IR_{\an,\exp}$ and therefore o-minimal. 
This shows again that the class of log-analytic functions is a proper extension of the class of globally subanalytic functions but shares its properties from the viewpoint of analysis and of o-minimality.
The latter observation does not hold when the exponential function comes into the game as remarked at the end of [4].
We formulate it in the following way:

\vs{0.5cm}
{\bf Example}

\vs{0.1cm}
Consider the function 
$$f: \IR\times\IR \to \IR,
(x,y)\mapsto 
\left\{\begin{array}{lll}| {y} |^{|2x|},&& y \neq 0, \\
&\mbox{if}&\\
0,&& y=0.
\end{array}\right.$$
which is definable in $\IR_\exp$.
Then the set of all $x\in \IR$ such that $f(x,-)$ is real-analytic is the set of integers. 

\vs{0.5cm}
The paper is organized as follows.
After a preliminary section on o-minimality, notations and conventions we give in Section 1 the definition of log-analytic functions and formulate elementary properties.
In Section 2 we present the preparation theorem of Lion-Rolin for log-analytic functions and prove the pure preparation theorem on simple sets.
Section 3 is devoted to the proof of the above theorems.

\section*{Preliminaries}

\subsection*{O-Minimality}

We give the definition and examples of o-minimal structures.

\vs{0.3cm}
{\bf Semialgebraic sets:}

A subset $A$ of $\IR^n, n\geq 1,$ is called {\bf semialgebraic} if there are $k,l\in \IN_0$ and real polynomials $f_i, g_{i1},\ldots,$ $ g_{ik}\in \IR[X_1,\ldots,X_n]$ for $1\leq i\leq l$, such that 
$$A=\bigcup_{i=1}^l \big\{x\in \IR^n\;\big\vert\; f_i(x)=0, g_{i,1}(x)>0,\ldots,g_{i,k}(x)>0\big\}.$$
A map is called semiagebraic if its graph is semialgebraic. 

\vs{0.3cm}
{\bf Semi- and subanalytic sets:}

A subset $A$ of $\IR^n, n\geq 1$, is called {\bf semianalytic} if for each $a\in \IR^n$ there are open neighbourhoods $U,V$ of $a$ with $\overline{U}\subset V$, $k,l\in \IN_0$ and real analytic functions $f_i,g_{i1},\ldots,g_{ik}$ on $V$ for $1\leq i\leq l$, such that
$$A\cap U=\bigcup_{i=1}^l\big\{x\in U\;\big\vert\; f_i(x)=0, g_{i,1}(x)>0,\ldots,g_{i,k}(x)>0\big\}.$$
A subset $B$ of $\IR^n, n\geq 1,$ is called {\bf subanalytic} if for each $a\in \IR^n$ there is an open neighbourhood $U$ of $a$, some $p\geq n$ and some bounded semianalytic set $A\subset\IR^p$ such that $B\cap U=\pi_n(A)$ where $\pi_n:\IR^p\to \IR^n,(x_1,\ldots,x_p)\mapsto (x_1,\ldots,x_n),$ is the projection on the first $n$ coordinates.

A map is called semianalytic or subanalytic if its graph is a semianalytic or subanalytic set, respectively. A set is called {\bf globally semianalytic} or {\bf globally subanalytic} if it is semianalytic or subanalytic, respectively, in the ambient projective space (or equivalently, after applying the semialgebraic homeomorphism $\IR^n\to \;]-1,1[^n, x_i\mapsto x_i/\sqrt{1+x_i^2}$).

\vs{0.3cm}
{\bf O-minimal structures:}

A {\bf structure} on the ordered field of reals is axiomatically defined as follows.

For $n\in \IN$ let $M_n$ be a set of subsets of $\IR^n$ and let $\mathcal{M}:=(M_n)_{n\in \IN}$. Then $\mathcal{M}$ is a structure on $\IR$ if the following holds for all $m,n,p\in \IN$.
\begin{itemize}
	\item[(S1)] If $A,B\in \IR^n$ then $A\cup B,A\cap B, \IR^n\setminus A\in M_n$. 
	\item[(S2)] If $A\in M_n$ and $B\in M_m$ then $A\times B\in M_{n+m}$.
	\item[(S3)] If $A\in M_p$ and $p\geq n$ then $\pi_n(A)\in M_n$.
	\item[(S4)] $M_n$ contains the semialgebraic subsets of $\IR^n$.
\end{itemize}
The structure $\mathcal{M}=(M_n)_{n\in\IN}$ on $\IR$ is called {\bf o-minimal} if additionally the following holds.
\begin{itemize}
	\item[(O)] The sets in $M_1$ are exactly the finite unions of intervals and points.
\end{itemize}
A subset of $\IR^n$ is called {\bf definable} in the structure $\mathcal{M}$ if it belongs to $M_n$. A function is definable in $\mathcal{M}$ if its graph is definable in $\mathcal{M}$.

The o-minimality axiom (O) implies that a subset of $\IR$ which is definable in an o-minimal structure on $\IR$ has only finitely many connected components. But much more can be deduced from o-minimality. A definable subset of 
 $\IR^n$, $n\in \IN$ arbitrary, has only finitely many connected components and these are again definable. More generally, sets and functions definable in an o-minimal structure exhibit tame geometric behaviour, for example the existence of definable cell decomposition. We refer to the book of Van den Dries [2] or to Miller and Van den Dries [5] for this and more of the general properties of o-minimal structures.

\vs{0.3cm}
{\bf Examples of o-minimal structures}
\begin{itemize}
	\item[(1)] The smallest o-minimal structure on $\IR$ is given by the semialgebraic sets. It is denoted by $\IR$.
	\item[(2)] $\IR_\exp$, the structure generated on  the real field by the global exponential function $\exp:\IR\to \IR_{>0}$ (i.e. the smallest structure containing the semialgebraic sets and the graph of the exponential function), is o-minimal. 
	\item[(3)] $\IR_\an$, the structure generated on the real field by the restricted analytic functions, is o-minimal. A function $f:\IR^n\to \IR$ is called restricted analytic if there is a function $g$ that is real analytic on a neighbourhood of $[-1,1]^n$ such that $f=g$ on $[-1,1]^n$ and $f=0$ else.
	The sets definable in $\IR_\an$ are exactly the globally subanalytic sets. 
	\item[(4)] $\IR_{\an,\exp}$, the structure generated by $\IR_\an$ and $\IR_\exp$, is o-minimal. 
\end{itemize}

\newpage
\subsection*{\bf Notations}

\vs{0.1cm}
The empty sum is by definition $0$ and the empty product is by definition $1$.

By $\IN=\{1,2,\ldots\}$ we denote the set of natural numbers and by $\IN_0=\{0,1,2,\ldots\}$ the set of nonnegative integers.
Given $x\in \IR$ let $\lceil x\rceil$ be the smallest integer which is not smaller than $x$. 

We set $\IR_{>0}:=\{x\in \IR\mid x>0\}$. For $m,n\in \IN$ we denote by $M(m\times n,\IR)$ the set of $m\times n$-matrices with real entries.
For $P\in M(m\times n,\IR)$ we denote by $^tP\in M(n\times m,\IR)$ its transpose.
Given $x\in \IR\setminus\{0\}$ we denote by $\mathrm{sign}(x)\in \{\pm 1\}$ its sign. 
For $a,b\in \IR$ with $a\leq b$ we denote by $[a,b]$ the closed interval and by $]a,b[$ the open interval with endpoints $a,b$, respectively.
Denoting by $|\;|$ the euclidean norm on $\IR^n$ we set $S^{n-1}:=\{x\in \IR^n\mid |x|=1\}$ and, for $a\in \IR^n$ and $\varepsilon\in \IR_{>0}$,
$B(a,\varepsilon):=\{x\in \IR^n\mid |x-a|<\varepsilon\}$. 
Given a subset $A$ of $\IR^n$ we denote by $\overline{A}$ its closure.
A domain in $\IR^n$ is a nonempty, open and connected subset of $\IR^n$.

We use the usual $o$-notation and $O$-notation. By the symbol $\sim$ we denote asymptotic eqivalence.
By $\log_k$ we denote the $k$-times iterated of the logarithm.

\subsection*{Conventions}

From now on definable means definable in the o-minimal structure $\IR_{\an,\exp}$.
A definable cell (see [2, 5]) is assumed to be analytic.

\section{Log-analytic functions}

We give the definition of a logarithmic-analytic function (log-analytic for short) established by Lion and Rolin in [8].

\vs{0.5cm}
{\bf 1.1 Definition} 

\vs{0.1cm}
	Let $X \subset \IR^m$ be definable and let $f:X \to \IR$ be a function. 
	\begin{itemize}
		\item[(a)] Let $r \in \IN_0$. By induction on $r$ we define that $f$ is {\bf log-analytic of order at most $r$}.
	
	\vs{0.2cm}
	{\bf Base case:} 
		The function $f$ is log-analytic of order at most $0$ if $f$ is piecewise the restriction of globally subanalytic functions; i.e. there is a finite decomposition $\ma{C}$ of $X$ into definable sets such that for $C \in \ma{C}$ there is a globally subanalytic function $F:\IR^m \to \IR$ such that $f|_C=F|_C$.
	
	\vs{0.2cm}
	{\bf Inductive step:}
		The function $f$ is log-analytic of order at most $r$ if the following holds. There is a finite decomposition $\ma{C}$ of $X$ into definable sets such that for $C \in \mathcal{C}$ there are $k,l \in \IN_0$, a globally subanalytic function $F:\IR^{k+l} \to \IR$ and log-analytic functions $g_1,...,g_k: C \to \IR, h_1,\ldots,h_l:C\to \IR_{>0}$ of order at most $r-1$ such that
		$$f|_C=F\big(g_1,...,g_k,\log(h_1),...,\log(h_l)\big).$$
	\item[(b)] Let $r\in \IN_0$. We call $f$ {\bf log-analytic of order $r$} if $f$ is log-analytic of order at most $r$ but not of order at most $r-1$.
	\item[(c)] We call $f$ {\bf log-analytic} if it is log-analytic of order $r$ for some $r\in \IN_0$.  
	\end{itemize}

\vs{0.2cm}
{\bf 1.2 Remark}

\begin{itemize}
	\item[(1)] A log-analytic function is definable.
	\item[(2)] Let $\ma{L}_\an$ be the natural language for $\IR_\an$; i.e. $\ma{L}_\an$ is the augmentation of the language $\ma{L}$ of ordered rings by symbols for restricted analytic functions. Then the log-analytic functions are precisely those definable functions which are piecewise given by $\ma{L}_\an(^{-1},(\sqrt[n]{...})_{n=2,3,...},\log)$-terms (compare with Van den Dries et al. [3]).
	\item[(3)] A function is log-analytic of order $0$ if and only if it is piecewise the restriction of globally subanalytic functions.
\end{itemize} 

\vs{0.2cm}
{\bf 1.3 Remark}

\vs{0.1cm}
Let $X\subset \IR^n$ be definable.
\begin{itemize}
	\item[(1)] Let $r\in \IN_0$. The set of log-analytic functions on $X$ of order at most $r$ is a ring with respect to pointwise addition and multiplication.
	\item[(2)] The set of log-analytic functions on $X$ is a ring with respect to pointwise addition and multiplication. 
\end{itemize}

\section{Preparation of log-analytic functions}

We let $x=(x_1,\ldots,x_n)$ range over $\IR^n$ and $y$ over $\IR$. 
We set $\pi:\IR^n\times \IR\to \IR^n, (x,y)\mapsto x$.

\subsection{Lion-Rolin preparation}

Lion and Rolin [8] have established a preparation result for log-analytic functions which we state here (in the corrected form of Paw\l ucki and Pi\c{e}kosz [11]).

Let $r\in \IN_0$. We let $w=(w_0,\ldots,w_r)$ range over $\IR^{r+1}$. 
We set
$\pi^*:\IR^n\times\IR^{r+1}\to \IR^n, (x,w)\mapsto x$.

\vs{0.5cm}
{\bf 2.1 Definition}

\vs{0.1cm}
Let $r\in \IN_0$ and let $D\subset \IR^n\times \IR^{r+1}$ be definable.  A function $u:D \to \IR$ is called a {\bf special unit} on $D$ if $u=v \circ \varphi$ where the following holds:
	\begin{itemize}
		\item [(a)] The function $\varphi$ is given by
		\begin{eqnarray*}
		\varphi:D &\to& [-1,1]^s,  \\
	 \varphi(x,w)	&=&\Big(b_1(x)\prod_{j=0}^r\vert{w_j}\vert^{p_{1j}},\ldots, b_s(x)\prod_{j=0}^r\vert{w_j}\vert^{p_{sj}}\Big),
	\end{eqnarray*}
		where $s\in \IN_0$, $p_{ij} \in \IQ$ for $(i,j) \in \{1,...,s\} \times \{0,...,r\}$ and $b_1,...,b_s$ are definable functions on $\pi^*(D)$ which have no zeros.
		\item [(b)] The function $v$  is a real power series which converges absolutely on an open neighbourhood of $[-1,1]^s$.
		\item [(c)] There are $d_1,d_2\in \IR_{>0}$ such that 
		$d_1\leq v\leq d_2$ on $[-1,1]^s$.
	\end{itemize} 
	We call $b:=(b_1,...,b_s)$ a tuple of {\bf base functions} for $u$ and 
	$\ma{I}:=\big(s,v,b,P\big)$
	where 
	$$P:=\left(\begin{array}{cccc}
	p_{10}&\cdot&\cdot&p_{1r}\\
	\cdot&& &\cdot\\
	\cdot&& &\cdot\\
    p_{s0}&\cdot&\cdot&p_{sr}\\
    \end{array}\right)\in M\big(s\times (r+1),\IR\big)$$
	a {\bf describing tuple} for $u$. 

\vs{0.5cm}
Let $C\subset \IR^n\times\IR$ be definable and let $r\in \IN_0$.

\vs{0.5cm}
{\bf 2.2 Definition}

\vs{0.1cm}
A tuple $\ma{Y}=(y_0,...,y_r)$ of functions on $C$ is called an {\bf $r$-logarithmic scale} on $C$ with {\bf center} $\Theta=(\Theta_0,\ldots,\Theta_r)$ if the following holds:
	\begin{itemize}
		\item[(a)] For every $j\in \{0,\ldots,r\}$ we have $y_j>0$ or $y_j<0$.
		\item[(b)] The functions $\Theta_j$ are definable functions on $\pi(C)$ for every $j\in \{0,\ldots,r\}$. 
	    \item[(c)] It is $y_0(x,y)=y-\Theta_0(x)$ and $y_j(x,y)=\log(|y_{j-1}(x,y)|)-\Theta_j(x)$ for $j\in \{1,\ldots,r\}$.
		\item[(d)] There is $\varepsilon_0 \in \,]0,1[$ such that $0<|y_0(x,y)|<\varepsilon_0|y|$ for all $(x,y) \in C$ or $\Theta_0 = 0$, and for every $j \in \{1,...,r\}$ there is $\varepsilon_j \in\, ]0,1[$ such that $0<|y_j(x,y)|<\varepsilon_j|\log(|y_{j-1}(x,y)|)|$ for all $(x,y) \in C$ or $\Theta_j = 0$.
	\end{itemize}

\vs{0.2cm}
Note that $\Theta$ is uniquely determined by $\ma{Y}$. Note also that $\ma{Y}$ is log-analytic if and only if $\Theta$ is log-analytic.





\vs{0.5cm}
{\bf 2.3 Definition}

\vs{0.1cm}
Let $\ma{Y}=(y_0,\ldots,y_r)$ be an $r$-logarithmic scale on $C$. Its {\bf sign} $\mathrm{sign}(\ma{Y})\in \{-1,1\}^{r+1}$ is defined by
$$\mathrm{sign}(\ma{Y})=\big(\mathrm{sign}(y_0),\ldots,\mathrm{sign}(y_r)\big).$$

\vs{0.5cm}
{\bf 2.4 Definition}

\vs{0.1cm}
Let $\ma{Y}=(y_0,\ldots,y_r)$ be an $r$-logarithmic scale on $C$.
Let $q=(q_0,\ldots,q_r)\in \IQ^{r+1}$.
We set
$$|\ma{Y}|^{\otimes q}:=\prod_{j=0}^r|y_j|^{q_j}.$$

\vs{0.5cm}
{\bf 2.5 Definition}

\vs{0.1cm}
Let $\ma{Y}$ be an $r$-logarithmic scale on $C$.
We set
$$C^\ma{Y}:= \big\{\big(x,\ma{Y}(x,y)\big) \;\big\vert\;(x,y) \in C\big\}\subset \IR^n\times\IR^{r+1}.$$
		
\vs{0.5cm}
{\bf 2.6 Definition}

\vs{0.1cm}
Let $g:C\to \IR$ be a function.
We say that $g$ is {\bf $r$-log-analytically prepared with respect to $y$} if 
$$g(x,y)=a(x)|\ma{Y}(x,y)|^{\otimes q}  u\big(x,\ma{Y}(x,y)\big)$$
where $a$ is a definable function on $\pi(C)$ which vanishes identically or has no zero, $\ma{Y}$ is an $r$-logarithmic scale on $C$, $q\in\IQ^{r+1}$ and $u$ is a special unit on $C^\ma{Y}$.
The  function $a$ is called {\bf coefficient} of $g$ and base functions $b_1,...,b_s$ of $u$ are called {\bf base functions} of $g$. 
We call
$\ma{J}:=\big(r,\ma{Y},a,q,\ma{I}\big)$ 
where $\ma{I}$ is a describing tuple for $u$ an {\bf $r$-preparing tuple} for $g$.

\vs{0.5cm}
{\bf 2.7 Remark}

\vs{0.1cm}
Let $g:C\to \IR$ be a function.
If $g$ is $r$-log-analytically prepared with respect to $y$ then $g$ is definable but not necessarily log-analytic.

\vs{0.5cm}
{\bf 2.8 Fact} [8, Section 2.1]

\vs{0.1cm}
{\it Let $X \subset \IR^n \times \IR$ be definable and let $f:X \to \IR$ be a log-analytic function of order $r$. 
	Then there is a definable cell decomposition $\ma{C}$ of $X$ such that $f|_C$ is $r$-log-analytically prepared with respect to $y$ for every $C \in \mathcal{C}$.}

\vs{0.5cm}
Here is the promised  example that the above can in general not be carried out in the log-analytic category.

\vs{0.5cm}
{\bf 2.9 Example}

\vs{0.1cm}
Let $\varphi:\,]0,+\infty[\,\to \IR, t\mapsto t/(1+t)$. 
Consider the log-analytic function 
$$f:\IR_{>0}\times \IR, (x,y)\mapsto -\frac{1}{\log(\varphi(y))}-x.$$ 
Then $f$ is log-analytic of order $1$ but does not allow a piecewiese $1$-log-analytic preparation with log-analytic data only.

\vs{0.1cm}
{\bf Proof:}

\vs{0.1cm}
Assume that the contrary holds.
Let $\ma{C}$ be a corresponding cell decomposition.
Set $\psi:\,]0,1[\to\IR, t\mapsto t/(1-t)$. Then $\psi$ is the compositional inverse of $\varphi$.
Note that $f(x,\psi(e^{-1/x}))=0$ for all $x\in \IR_{>0}$.
Let $\alpha:\IR_{>0}\to \IR, x\mapsto \psi(e^{-1/x})$. Then $\alpha$ is not log-analytic and $\alpha(x)=\sum_{n=1}^\infty e^{-n/x}$ for all $x\in \IR_{>0}$.
By passing to a finer cell decomposition we find
a cell $C$ of the form
$$C:=\Big\{(x,y)\in \IR_{>0}\times\IR_{>0}\;\big\vert\; 0<x<\varepsilon, \alpha(x)<y<\alpha(x)+\eta(x)\Big\}$$
with some suitable $\varepsilon\in \IR_{>0}$ and some definable function $\eta:\,]0,\varepsilon[\to \IR_{>0}$ such that $f$ is $1$-log-analytically prepared on $C$ with log-analytic data only.
Let $\big(r,\ma{Y},a,q,\ma{I}\big)$ be a $1$-preparing tuple for $f|_C$ and let 
$\Theta=(\Theta_0,\Theta_1)$ be the center of $\ma{Y}$. 

\vs{0.2cm}
{\bf Claim:} $\Theta_0=0$.

\vs{0.1cm}
{\bf Proof of the claim:} Assume that $\Theta_0$ is not the zero function. By the definition of a $1$-logarithmic scale we find $\varepsilon_0\in\,]0,1[$ such that $|y_0|<\varepsilon_0 |y|$ on $C$. This implies $|\alpha(x)-\Theta_0(x)|\leq \alpha(x)$ for all $0<x<\varepsilon$. But this is not possible since we have $\alpha=o(\Theta_0)$ at $0$ by the assumption that $\Theta_0$ is log-analytic and not the zero function.
\hfill$\blacksquare_{\mathrm{Claim}}$

\vs{0.2cm} 
From 
$$f(x,y)=a(x)|\ma{Y}(x,y)|^{\otimes q}u(x,\ma{Y}(x,y))$$ 
for all $(x,y)\in C$ and 
$\lim_{y\searrow \alpha(x)}f(x,y)=0$ for all $x\in \,]0,\varepsilon[$ we get by o-minimality that
there is, after shrinking $\varepsilon>0$ if necessary, some $j\in \{0,1\}$ such that
$\lim_{y\searrow \alpha(x)}|y_j(x,y)|^{q_j}=0$ for all 
$x\in \,]0,\varepsilon[$. By the claim the case $j=0$ is not possible. 
In the case $j=1$ we have, again by the claim, that $q_1>0$ and therefore
$\Theta_1=\log(\alpha)$. But this is a contradiction to the assumption that the function $\Theta_1$ is log-analytic since the function $\log(\alpha)$ on the right is not log-analytic. This can be seen by applying the logarithmic series. We obtain that
$\log(\alpha(x))+1/x\sim e^{-1/x}$.
\hfill$\blacksquare$

\vs{0.5cm}
In the case of globally subanalytic functions (i.e. log-analytic functions of order $0$) there is the following well known stronger result.

\vs{0.5cm}
{\bf 2.10 Remark} [8, Section 1]

\vs{0.1cm}
If $f:X\to \IR$ is globally subanalytic then 
there is a globally subanalytic cell decomposition $\ma{C}$ of $X$ such that $f|_C$ is globally subanalytic prepared with respect to $y$ for all $C\in\ma{C}$; i.e.
$$f|_C(x,y)=a(x)|y-\Theta(x)|^q v\big(b_1(x)|y-\Theta(x)|^{p_1},\ldots,b_s(x)|y-\Theta(x)|^{p_s}\big)$$
where $a,b_1,\ldots,b_s,\Theta:\pi(C)\to \IR$ are globally subanalytic functions such that the image of $\big(b_1(x)|y-\Theta(x)|^{p_1},\ldots,b_s(x)|y-\Theta(x)|^{p_s}\big)$ is contained in $[-1,1]^s$ and $v$ is a real power series converging absolutely on $[-1,1]^s$ with $v([-1,1]^s)\subset \IR_{>0}$.

\subsection{Simple sets and simple preparation}

Let $C\subset \IR^n\times\IR$ be definable and $r\in \IN_0$.

\vs{0.5cm}
{\bf 2.11 Definition}

\begin{itemize}
\item[(a)]
We call $C$ {\bf $r$-admissible} if there is an $r$-logarithmic scale on $C$. 
\item[(b)] 
We call $C$ {\bf $r$-unique} if there is a unique $r$-logarithmic scale on $C$. 
\end{itemize}

\vs{0.2cm}
{\bf 2.12 Definition}

\vs{0.1cm}
An $r$-logarithmic scale on $C$ is called {\bf elementary} if its center is vanishing. 

\vs{0.5cm}
{\bf 2.13 Remark}

\vs{0.1cm}
An elementary $r$-logarithmic scale may not exist on $C$. If it exists it is uniquely determined and log-analytic.

\vs{0.5cm}
{\bf 2.14 Definition}

\vs{0.1cm}
If $C$ has an elementary $r$-logarithmic scale then we call $C$ {\bf $r$-elementary}. The elementary $r$-logarithmic scale on $C$ is then denoted by $\ma{Y}_r^\mathrm{el}=\ma{Y}_{r, C}^\mathrm{el}$.

\vs{0.5cm}
For the next definition compare with the setting of [4, Section 4] and [7, Definition 3.6]. Given $x\in \IR^n$, we set $C_x:=\{y\in \IR\mid (x,y)\in C\}$.

\vs{0.5cm}
{\bf 2.15 Definition}

\vs{0.1cm}
We call $C$ {\bf simple} if for every $x\in \pi(C)$ we have $C_x=\,]0,d_x[$ for some $d_x\in \IR_{>0}\cup\{+\infty\}$.

\vs{0.5cm}
{\bf 2.16 Remark}

\vs{0.1cm}
Let $\ma{C}$ be a definable cell decomposition of $\IR^n\times\IR_{>0}$. Then
$$\IR^n=\bigcup\{\pi(C)\mid C\in\ma{C}\mbox{ simple}\}.$$

\vs{0.5cm}
We set $e_0:=0$ and $e_r:=\exp(e_{r-1})$ for $r\in \IN$. In the following we let $1/0:=+\infty$.

\vs{0.5cm}
{\bf 2.17 Proposition}

\vs{0.1cm}
{\it Let $C$ be simple and $r$-elementary and let $\ma{Y}_{r,C}^\mathrm{el}=(y_0,\ldots,y_r)$. Then the following holds:
\begin{itemize}
	\item[(1)] $\sup C_x\leq 1/e_r$ for all $x\in \pi(C)$.
    \item[(2)] $y_0=y, y_1=\log(y), y_j=\log_{j-1}(-\log(y))$ for $j\in \{2,\ldots,r\}$.
	\item[(3)] $\mathrm{sign}(\ma{Y}_{r,C}^\mathrm{el})=(1,-1,1,\ldots,1)\in \IR^{r+1}$.
\end{itemize}}
{\bf Proof:}

\vs{0.1cm}
For $x\in \pi(C)$ let $d_x:=\sup C_x$.
Let  
$(\sigma_0,\ldots,\sigma_r)$ be the sign of $\ma{Y}_{r,C}^\mathrm{el}$.  
We show by induction on $k\in \{0,\ldots,r\}$ that $d_x\leq 1/e_k$ for all $x\in \pi(C)$, that $y_0=y,y_1=\log(y),y_j=\log_{j-1}(-\log(y))$ for all $j\in \{2,\ldots,k\}$ and that $(\sigma_0,\ldots,\sigma_k)=(1,-1,1,\ldots,1)\in \IR^{k+1}$.

\vs{0.2cm}
$k=0$:
We have $y_0=y$ by Definition 2.2. This gives $\sigma_0=1$. That $d_x\leq +\infty=1/e_0$ for all $x\in \pi(C)$ is clear.

\vs{0.2cm}
$k=1$:
Since $y_0=y$ and $\sigma_0=1$ by above we obtain according to Definition 2.2 that $y_1=\log(y)$ and that $d_x\leq 1=1/e_1$ for all $x\in \pi(C)$.
This gives $\sigma_1=-1$.

\vs{0.2cm}
$k=2$:
Since $\sigma_1=-1$ we have $y_1<0$. According to Definition 2.2  we get that $y_2=\log(-y_1)=\log(-\log(y))$ and therefore that $d_x\leq 1/\exp(1)=1/e_2$ and $\sigma_2=1$.

\vs{0.2cm}
$k\to k+1$:
We can assume that $k\geq 2$. By the inductive hypothesis we have $y_k=\log_{k-1}(-\log(y))>0$ and $\sigma_k=1$.
According to Definition 2.2 we obtain that $y_{k+1}=\log_k(-\log(y))$ and that $d_x\leq 1/\exp(e_k)=1/e_{k+1}$ for all $x\in \pi(C)$. 
This gives also $\sigma_{k+1}=1$.
\hfill$\blacksquare$

\vs{0.5cm}
{\bf 2.18 Definition}

\vs{0.1cm}
We call $C$ {\bf $r$-simple} if it is simple and $r$-admissible.

\vs{0.5cm}
{\bf 2.19 Proposition}

\vs{0.1cm}
{\it Let $C$ be $r$-simple. Then $C$ is $r$-elementary and $r$-unique.}

\vs{0.1cm}
{\bf Proof:}

\vs{0.1cm}
Let $\ma{Y}=(y_0,\ldots,y_r)$ be an $r$-logarithmic scale on $C$.  
We show that $\ma{Y}$ is elementary and are done by Remark 2.13.
Let  
$\Theta=(\Theta_0,\ldots,\Theta_r)$ be the center of $\ma{Y}$.  
We show  by induction on $k\in \{0,\ldots,r\}$
that $\Theta_0=\ldots=\Theta_k=0$.

\vs{0.2cm}
$k=0$:
Assume that $\Theta_0 \neq 0$. Then by Definition 2.2 there is $\varepsilon_0 \in\, ]0,1[$ such that
$$|y-\Theta_0(x)| < \varepsilon_0|y|$$
for all $(x,y) \in C$. Let $x \in \pi(C)$ such that $\Theta_0(x) \neq 0$. Then we obtain
$$+\infty=\lim_{y\searrow 0}\Big\vert 1-\frac{\Theta_0(x)}{y}\Big\vert\leq \varepsilon_0$$ 
which is a contradiction.

\vs{0.2cm}
$k=1$:
Assume that $\Theta_1 \neq 0$.  By the case $k=0$ and Proposition 2.17 we have $y_0=y$. According to Definition 2.2 there is $\varepsilon_1 \in\, ]0,1[$ such that
$$|\log(y) - \Theta_1(x)| < \varepsilon_1 |\log(y)|$$
for all $(x,y)\in C$. 
Therefore 
$$1=\lim_{y\searrow 0}\Big\vert 1-\frac{\Theta_1(x)}{\log(y)}\Big\vert \leq \varepsilon_1$$
for $x\in \pi(C)$ 
which is a contradiction. 

\vs{0.2cm}
$k=2$:
Assume that $\Theta_2 \neq 0$.  By the case $k=1$ and Proposition 2.17 we have $y_1=\log(y)$. According to Definition 2.2 there is $\varepsilon_2 \in\, ]0,1[$ such that
$$|\log(-y_1) - \Theta_1(x)| < \varepsilon_2 |\log(-y_1)|$$
for all $(x,y)\in C$. 
Therefore 
$$1=\lim_{y\searrow 0}\Big\vert 1-\frac{\Theta_1(x)}{\log(-y_1)}\Big\vert \leq \varepsilon_2$$
for $x\in \pi(C)$ 
which is a contradiction.

\vs{0.2cm}
$k\to k+1$:
We may assume that $k\geq 2$.
Assume that $\Theta_{k+1} \neq 0$. By the inductive hypothesis and Proposition 2.17 we have $y_k=\log_{k-1}(-\log(y))$. 
According to Definition 2.2 there is $\varepsilon_{k+1} \in\, ]0,1[$ such that
$$|\log(y_k)-\Theta_{k+1}(x)| < \varepsilon_{k+1}|\log(y_k)|$$ 
for all $(x,y) \in C$. Therefore 
$$1=\lim_{y \searrow 0}\Big\vert 1-\frac{\Theta_{k+1}(x)}{\log(y_k)}\Big\vert< \varepsilon_{k+1}$$
for $x\in \pi(C)$ which is a contradiction.  
\hfill$\blacksquare$

\vs{0.5cm}
{\bf 2.20 Corollary}

\vs{0.1cm}
{\it Let $C$ be simple. Then the following are equivalent:
	\begin{itemize}
		\item[(i)] $C$ is $r$-simple.
		\item[(ii)] $\sup C_x\leq 1/e_r$ for every $x\in \pi(C)$.
\end{itemize}}
{\bf Proof:}

\vs{0.1cm}
(i) $\Rightarrow$ (ii): If $C$ is $r$-simple then $C$ is $r$-elementary by Proposition 2.19. By Proposition 2.17 we obtain (ii).

\vs{0.2cm}
(ii) $\Rightarrow$ (i):
Let $y_0=y, y_1:=\log(y), y_j=\log_{j-1}(-\log(y))$ for $j\in \{2,\ldots,r\}$. Then it is straightforward to see that $\ma{Y}:=(y_0,\ldots,y_r)$ is a well defined elementary $r$-logarithmic scale on $C$.
\hfill$\blacksquare$

\vs{0.5cm}
{\bf 2.21 Definition}

\vs{0.1cm}
Let $g:C\to \IR$ be a function. We say that $g$ is {\bf elementary $r$-log-analytically prepared with respect to $y$} if $g$ is $r$-log-analytically prepared with elementary $r$-logarithmic scale.

\vs{0.5cm}
{\bf 2.22 Corollary}

\vs{0.1cm}
{\it Let $C$ be simple and let $g:C\to \IR$ be $r$-log-analytically prepared with respect to $y$. Then $g$ is elementary $r$-log-analytically prepared with respect to $y$.}

\vs{0.1cm}
{\bf Proof:}

\vs{0.1cm}
Since $C$ is simple and since $g$ has an $r$-log-analytic preparation 
we have that $C$ is $r$-simple. We are done by Proposition 2.19.
\hfill$\blacksquare$

\vs{0.5cm}
Let $q=(q_0,\ldots,q_r)\in \IQ^{r+1}$ with $q\neq 0$.
We set
$j(q):=\min\{j\mid q_j\neq 0\}$ and $\sigma(q):=\mathrm{sign}(q_{j(q)})\in \{\pm 1\}$. 
Moreover, let
$$q_\mathrm{diff}:=\big(q_0-1,\ldots,q_{j(q)}-1,q_{j(q)+1},\ldots,q_r\big).$$

\vs{0.5cm}
{\bf 2.23 Remark}

\vs{0.1cm}
Let $C$ be $r$-simple. Let 
	$q\in \IQ^{r+1}$ with $q\neq 0$.
	Then
	$$\lim_{y\searrow 0}|\ma{Y}_{r,C}^\mathrm{el}(y)|^{\otimes q}=
	\left\{\begin{array}{cc}
	0,&j(q)=0,\sigma(q)=+1,\\
	+\infty,&j(q)=0,\sigma(q)=-1,\\
	+\infty,&j(q)>0,\sigma(q)=+1,\\
	0,&j(q)>0,\sigma(q)=-1.\\
	\end{array}\right.$$

\vs{0.5cm}
{\bf 2.24 Proposition}

\vs{0.1cm}
{\it Let $C$ be $r$-simple. Let $q\in \IQ^{r+1}$ with $q\neq 0$. 
Then
$$\lim_{y \searrow 0}\Big\vert\frac{\frac{d}{dy}|\ma{Y}_{r,C}^\mathrm{el}(y)|^{\otimes q}}
{|\ma{Y}_{r,C}^\mathrm{el}(y)|^{\otimes q_\mathrm{diff}}}\Big\vert=q_{j(q)}.$$}

\vs{0.1cm}
{\bf Proof:}

\vs{0.1cm}
Let $\ma{Y}_{r,C}^\mathrm{el}=(y_0,\ldots,y_r)$.
We get by Proposition 2.17 that $d|y_0|/dy=1$ and that for $j\in \{1,\ldots,r\}$
$$\frac{d|y_j|}{dy}=-\frac{1}{\prod_{i=0}^{j-1}|y_i|}.$$
This gives
\begin{eqnarray*}
\frac{d}{dy}\prod_{j=0}^r|y_j|^{q_j}&=&q_0|y_0|^{q_0-1}|y_1|^{q_1}\cdot\ldots\cdot |y_r|^{q_r}-\sum_{j=1}^rq_j|y_j|^{q_j-1}\frac{1}{\prod_{i=0}^{j-1}|y_i|}\prod_{i\neq j}|y_i|^{q_i}\\
&=&q_0|y_0|^{q_0-1}\prod_{i>0}|y_i|^{q_i}-\sum_{j=1}^rq_j\prod_{i\leq j}|y_i|^{q_i-1}\prod_{i>j}|y_i|^{q_i}.\\
\end{eqnarray*}	
We obtain the assertion by the growth properties of the iterated logarithms.

\hfill$\blacksquare$

\subsection{Pure preparation}

Let $C\subset \IR^n\times\IR$ be definable.	

\vs{0.5cm}
{\bf 2.25 Definition}

\vs{0.1cm}
Let $g:C\to \IR$ be a function.
We say that $g$ is {\bf purely $r$-log-analytically prepared with respect to $y$} if $g$ is $r$-log-analytically prepared with respect to $y$  
with log-analytic logarithmic scale, log-analytic coefficient and log-analytic base functions.
An $r$-preparing tuple for $g$ with log-analytic components is then called a {\bf purely $r$-preparing tuple} for $g$. 

\vs{0.5cm}
{\bf 2.26 Remark}

\vs{0.1cm}
Let $g:C\to \IR$ be a function. 
If $g$ is purely $r$-log-analytically prepared  with respect to $y$ then $g$ is log-analytic. 

\vs{0.5cm}
In the next proposition and corollary we assume that $r\geq 1$.

\vs{0.5cm}
{\bf 2.27 Proposition}

\vs{0.1cm}
{\it Let $C$ be $r$-simple.
\begin{itemize}
	\item[(1)] Let $g:C\to \IR$ be purely $(r-1)$-log-analytically prepared with respect to $y$.
	Then there are $t\in \IN$, a log-analytic function $\eta:\pi(C)\to \IR^t$ and a globally subanalytic function $G:\IR^t\times\IR^r\to \IR$
	such that
	$$g(x,y)=G(\eta(x),\ma{Y}_{r-1,C}^\mathrm{el}(y))$$
	for all $(x,y)\in C$.
	\item[(2)] Let $h:C\to \IR_{>0}$ be purely $(r-1)$-log-analytically prepared with respect to $y$. 
	Then there are $t\in \IN$, a log-analytic function $\eta:\pi(C)\to \IR^t$ and a globally subanalytic function $H:\IR^t\times \IR^{r+1}\to\IR$
	such that
	$$\log(h(x,y))=H(\eta(x),\ma{Y}_{r,C}^\mathrm{el}(y))$$ 
	for all $(x,y)\in C$.
\end{itemize}}
{\bf Proof:}

\vs{0.1cm}
(1): By Corollary 2.20 we have that $C$ is $(r-1)$-simple. By Corollary 2.22 we have
a purely $(r-1)$-preparing tuple $\ma{J}$ for $g$ of the form
$$\ma{J}=\big(r-1,\ma{Y}_{r-1,C}^\mathrm{el},a,q,s,v,b,P\big)$$ where $a,b_1,\ldots,b_s$ are log-analytic functions on $\pi(C)$. 
Take $t=s+1$ and  
$$\eta=(\eta_0,\ldots,\eta_s):\pi(C)\to \IR^t, x\mapsto \big(a(x),b_1(x),\ldots,b_s(x)\big).$$
Then $\eta$ is log-analytic.
Let $z=(z_0,\ldots,z_s)$ and $w=(w_0,\ldots,w_{r-1})$. 
Set
$$\alpha_0:\IR^t\times\IR^r\to \IR, (z,w)\mapsto z_0\prod_{j=0}^{r-1}|w_j|^{q_j}.$$
For $i\in \{1,\ldots,s\}$ let
$$\alpha_i:\IR^t\times\IR^r\to \IR, (z,w)\mapsto z_i\prod_{j=0}^{r-1}|w_j|^{p_{ij}}.$$
Set
$$G:\IR^t\times\IR^r\to \IR, (z_0,\ldots,z_s,w_0,\ldots,w_{r-1})\mapsto$$
$$\left\{\begin{array}{ccc}
\alpha_0(z,w)v\Big(\alpha_1(z,w),\ldots,\alpha_s(z,w)\Big),&|\alpha_i(z,w)|\leq 1\mbox{ for all }i\in \{1,\ldots,s\},\\
0,&\mbox{else}.\\
\end{array}\right.$$
Then $G$ is globally subanalytic and we have
$$g(x,y)=G(\eta(x),\ma{Y}_{r-1,C}^\mathrm{el}(y))$$
for all $(x,y)\in C$.

\vs{0.2cm}
(2): By Corollary 2.20 we have that $C$ is $(r-1)$-simple. By Corollary 2.22 we have
a purely $(r-1)$-preparing tuple $\ma{J}$ for $h$ of the form
$$\ma{J}=\big(r-1,\ma{Y}_{r-1,C}^\mathrm{el},a,q,s,v,b,P\big)$$ 
where $a,b_1,\ldots,b_s$ are log-analytic functions on $\pi(C)$. 
Then $a>0$.
Take $t=s+1$ and
$$\eta=(\eta_0,\ldots,\eta_s):\pi(C)\to \IR^t, x\mapsto \big(\log(a(x)),b_1(x),\ldots,b_s(x)\big).$$
Then $\eta$ is log-analytic.
Let $z=(z_0,\ldots,z_s)$ and $w=(w_0,\ldots,w_r)$. 
Set
$$\beta_0:\IR^t\times\IR^r\to \IR, (z,w)\mapsto z_0+\sum_{j=0}^{r-1}q_jw_{j+1}.$$
For $i\in \{1,\ldots,s\}$ let
$$\alpha_i:\IR^t\times\IR^r\to \IR, (z,w)\mapsto z_i\prod_{j=0}^{r-1}|w_j|^{p_{ij}}.$$
Set
$$H:\IR^t\times\IR^r\to \IR, (z_0,\ldots,z_s,w_0,\ldots,w_{r-1})\mapsto$$
$$\left\{\begin{array}{ccc}
\beta_0(z,w)+\log\Big(v\Big(\alpha_1(z,w),\ldots,\alpha_s(z,w)\Big)\Big),&|\alpha_i(z,w)|\leq 1\mbox{ for all }i\in \{1,\ldots,s\},\\

0,&\mbox{else.}\\
\end{array}\right.$$
Then $H$ is globally subanalytic since $\log(v)$ is globally subanalytic and we have
$$\log(h(x,y))=H(\eta(x),\ma{Y}_{r,C}^\mathrm{el}(y))$$
for all $(x,y)\in C$.
\hfill$\blacksquare$

\vs{0.5cm}
{\bf 2.28 Corollary} 

\vs{0.1cm}
{\it Let $C$ be $r$-simple and let $g_1,\ldots,g_k:C\to \IR$ and $h_1,\ldots,h_l:C\to \IR_{>0}$ be purely $(r-1)$-log-analytically prepared  with respect to $y$.
	Let $F:\IR^{k+l}\to\IR$ be globally subanalytic.
	Then there are $t\in \IN$, a log-analytic function $\eta:\pi(C)\to \IR^t$ and a globally subanalytic function $I:\IR^t\times \IR^{r+1}\to\IR$ such that
	$$F\big(g_1(x,y),\ldots,g_k(x,y),\log(h_1(x,y)),\ldots,\log(h_l(x,y))\big)=I(\eta(x),\ma{Y}_{r,C}^\mathrm{el}(x,y))$$
    for all $(x,y)\in C$.}

\vs{0.5cm}
{\bf 2.29 Proposition}

\vs{0.1cm}
{\it Let $C$ be $r$-simple and let $i:C\to \IR$ be a function. Assume that there are $t\in \IN$, a log-analytic function $\eta:\pi(C)\to \IR^t$ and a globally subanalytic function $I:\IR^t\times \IR^{r+1}\to \IR$ such that
$$i(x,y)=I(\eta(x),\ma{Y}_{r,C}^\mathrm{el}(y))$$
for all $(x,y)\in C$. Then there is a definable cell decomposition $\ma{D}$ of $C$ such that $i|_D$ is purely $r$-log-analytically prepared with respect to $y$ for every simple $D\in \ma{D}$.}

\vs{0.1cm}
{\bf Proof:}

\vs{0.1cm}
We do induction on $r$.
Let $z:=(z_1,\ldots,z_t)$ and $w:=(w_0,\ldots,w_r)$. We set $w':=(w_1,\ldots,w_r)$.

\vs{0.2cm}
$r=0$: 
We prepare $I$ globally subanalytically with respect to the coordinate $w_0$ and find a globally subanalytic cell decomposition $\ma{E}$ of $\IR^t\times \IR_{>0}$ such that for every 
$E\in \ma{E}$ the restriction $I|_E$ is globally subanalytic prepared with respect to $w_0$. 
We find a definable cell decomposition $\ma{A}$ of $C$ such that for every simple $A\in \ma{A}$ there is $E_A\in \ma{E}$ such that 
$(\eta(x),y)\in E_A$ for every $(x,y)\in A$.
Fix a simple $A\in \ma{A}$. 
Since $A$ is simple and $y$ is plugged in for $w_0$ we get that $E_A$ is simple with respect to $w_0$ and therefore by Remark 2.10 in combination with Corollary 2.22 that  $I|_{E_A}$ is elementarily globally subanalytic prepared with respect to $w_0$. 
Hence again by Remark 2.10 we have that
$$I|_{E_A}=a(z)|w_0|^q v\big(b_1(z)|w_0|^{p_1},\ldots,b_s(z)|w_0|^{p_s}\big)$$
where $a,b_1,\ldots,b_s$ are globally subanalytic functions on $\pi(E_A)$. 
We obtain
$$i|_C(x,y)=a(\eta(x))|y|^q v\big(b_1(\eta(x))|y|^{p_1},\ldots,b_s(\eta(x))|y|^{p_s}\big)$$
for all $(x,y)\in C$ and are done.

\vs{0.2cm}
$r-1\to r$:
We prepare $I$ globally subanalytically with respect to the coordinate $w_0$ and
find a globally subanalytic cell decomposition $\ma{E}$ of $\IR^t\times \IR_{>0}\times \IR^r$ such that for every 
$E\in \ma{E}$ the restriction $I|_E$ is globally subanalytic prepared with respect to $w_0$. 
We find a definable cell decomposition $\ma{A}$ of $C$ such that for every simple $A\in \ma{A}$ there is $E_A\in \ma{E}$ such that 
$(\eta(x),\ma{Y}_{r,C}^\mathrm{el}(y))\in E_A$ for every $(x,y)\in A$.
Fix a simple $A\in \ma{A}$. 
Since $A$ is simple and $y$ is plugged in for $w_0$ we get that $E_A$ is simple with respect to $w_0$ and therefore by Remark 2.10 in combination with Corollary 2.22 that  $I|_{E_A}$ is elementarily globally subanalytic prepared with respect to $w_0$. 
Hence again by Remark 2.10 we have that
$$I|_{E_A}=a(z,w')|w_0|^q v\big(b_1(z,w')|w_0|^{p_1},\ldots,b_s(z,w')|w_0|^{p_s}\big)$$
where $a,b_1,\ldots,b_s$ are globally subanalytic functions. Denoting by $c$ one of these functions we have to deal with $c(\eta(x),y_1,\ldots,y_r)$.
Using composition of power series we are done with the following claim.

\vs{0.2cm}
{\bf Claim:}
Let $J:\IR^t\times \IR^r\to \IR$ be globally subanalytic. Then
there is a definable cell decomposition $\ma{B}$ of $A$ such that for every 
simple $B\in \ma{B}$ the function
$$j:B\to \IR, (x,y)\mapsto J\big(\eta(x),y_1,\ldots,y_r\big)$$
is purely $r$-log-analytically prepared.

\vs{0.1cm}
{\bf Proof of the claim:}
Since $C$ is $r$-simple we get that $A$ is $r$-simple.
Set
$$\widehat{A}:=\big\{\big(x,-1/\log(y)\big)\;\big\vert\; (x,y)\in A\big\}.$$
Then $\widehat{A}$ is $(r-1)$-simple by Corollary 2.20.
Set
$$\widehat{J}:\IR^t\times \IR^r\to\IR,
(z,w')\mapsto
\left\{\begin{array}{ccc}
\alpha(z,-1/w_1,-w_2,w_3,\ldots,w_r),&&w_1<0,\\
&\mbox{if}&\\
0,&&\mbox{else}.
\end{array}\right.$$
Then for $(x,y)\in A$ we have
$$j(x,y)=\widehat{J}\Big(\eta(x),\ma{Y}_{r-1,\widehat{A}}^\mathrm{el}\big(-1/\log(y)\big)\Big).$$
Applying the inductive hypothesis to
$$\widehat{j}:\widehat{A}\to \IR, (x,y)\mapsto \widehat{J}\big(\eta(x),\ma{Y}_{r-1,\widehat{A}}^\mathrm{el}(y)\big)$$
we are done.
\hfill$\blacksquare_\mathrm{Claim}$

\hfill$\blacksquare$

\vs{0.5cm}
{\bf 2.30 Theorem}

\vs{0.1cm}
{\it Let $f:\IR^n\times\IR_{>0}\to \IR$ be log-analytic of order $r$. Then there is a definable cell decomposition $\ma{C}$ of $X$ such that for every simple $C\in \ma{C}$ the cell $C$ is $r$-simple and $f|_C$ is purely $r$-log-analytically prepared with respect to $y$.}

\vs{0.1cm}
{\bf Proof:}

\vs{0.1cm}
We do induction on the log-analytic order $r$ of $f$.

\vs{0.2cm}
$r=0$: Then $f$ is piecewise globally subanalytic and we are done by Remark 2.10.

\vs{0.2cm}
$<r\to r$: It is enough to consider the following situation. Let $g_1,\ldots,g_k:\IR^n\times\IR_{>0}\to\IR, h_1,\ldots,h_l:\IR^n\times\IR_{>0}\to \IR_{>0}$ be log-analytic functions of order at most $r-1$ and let 
$F:\IR^{k+l}\to\IR$ be globally subanalytic such that 
$$f=F(g_1,\ldots,g_k,\log(h_1),\ldots,\log(h_l)).$$
Applying the inductive hypothesis and Corollary 2.28 in combination with Corollary 2.20 we find a definable cell decomposition $\ma{C}$ of $\IR^n\times\IR_{>0}$ such that
every simple $C\in \ma{C}$ is $r$-simple and for every such $C$ there is $t\in \IN$, a log-analytic function $\eta:\pi(C)\to \IR^t$ of order at most $r$ and a globally subanalytic function 
$E:\IR^t\times\IR^{r+1}\to\IR, (z,w)\mapsto E(z,w),$ such that
$$f|_C=E(\eta(x),\ma{Y}_{r,C}^\mathrm{el}(x,y))$$ 
for all $(x,y)\in C$.
By Proposition 2.29 we can refine the cell decomposition such that the assertion follows. 
\hfill$\blacksquare$

\vs{0.5cm}
Note that the above result extends the Expansion Theorem of [4, Section 4] to log-analytic functions.

\newpage
\section{Proof of the main results}

\subsection{Derivatives - Proof of Theorem A}

\vs{0.5cm}
{\bf 3.1 Theorem}

\vs{0.1cm}
{\it Let $f:\IR\times \IR_{>0}\to \IR, (x,y)\mapsto f(x,y),$ be log-analytic. Assume that for every $x\in \IR^n$ we have $\lim_{y\searrow 0}f(x,y)\in\IR$.
	Then the function $h:\IR^n\to \IR, x\mapsto \lim_{y\searrow 0}f(x,y),$ is log-analytic.}

\vs{0.1cm}
{\bf Proof:}

\vs{0.1cm}
By Theorem 2.30 we find a definable cell decomposition $\ma{C}$ of $\IR^n\times \IR_{>0}$ such that for every simple $C\in\ma{C}$ the cell $C$ is $r$-simple and 
$f|_C$ is purely $r$-log-analytically prepared with respect to $y$.
Let $C\in \ma{C}$ be a simple cell. 
Set $g:=f|_C$ and let 
$$\ma{J}=\big(r,\ma{Y}_{r,C}^\mathrm{el},a,q,s,v,b,P\big)$$
be a purely $r$-preparing tuple for $g$. Then
$$g(x,y)=a(x)|\ma{Y}_{r,C}^{\mathrm{el}}(y)|^{\otimes q}v\Big(b_1(x)|\ma{Y}_{r,C}^\mathrm{el}(y)|^{\otimes p_1},\ldots,b_s(x)|\ma{Y}_{r,C}^\mathrm{el}(y)|^{\otimes p_s}\Big)$$
for $(x,y)\in C$. 
By the assumption, Remark 2.23 and Definition 2.1 we see that 
$$A:\pi(C)\to \IR, x\mapsto \lim_{y\searrow 0}a(x)|\ma{Y}_{r,C}^{\mathrm{el}}(y)|^{\otimes q},$$
and, for $j\in \{1,\ldots,s\}$,
$$B_i:\pi(C)\to [-1,1],x\mapsto \lim_{y\searrow 0}b_j(x)|\ma{Y}_{r,C}^{\mathrm{el}}(y)|^{\otimes p_i},$$
are well defined and log-analytic.
We obtain that for $x\in \pi(C)$
$$h(x)=A(x)v\Big(B_1(x),\ldots,B_s(x)\Big).$$
Hence $h$ is log-analytic on $\pi(C)$. 
By Remark 2.16 we get that 
$h$ is log-analytic.

\hfill$\blacksquare$

\vs{0.5cm}
With the above theorem we are able to establish Theorem A.

\vs{0.5cm}
{\bf Theorem A}

\vs{0.1cm}
{\it Let $U\subset \IR^n$ be definable and open and let $f:U\to \IR$ be log-analytic. Let $i\in \{1,\ldots,n\}$ be such that $f$ is differentiable with respect to the variable $x_i$ on $U$.  Then $\partial f/\partial x_i$ is log-analytic.}

\vs{0.1cm}
{\bf Proof:}

\vs{0.1cm}
We may assume that $f$ is differentiable with respect to the last variable $x_n$. We have to show that $\partial f/\partial x_n$ is log-analytic.
Let $\mathfrak{e}_n:=(0,\ldots,0,1)\in \IR^n$ be the $n^\mathrm{th}$ unit vector.
We define
$$V:=\big\{(x,y)\in U\times \IR_{>0}\;\big\vert\; x+y\mathfrak{e}_n\in U\}$$ and 
$$F:V\to \IR, (x,y)\mapsto 
\frac{f(x+y\mathfrak{e}_n)-f(x)}{y}.$$
Then $F$ is log-analytic. 
Since
$$\frac{\partial f}{\partial x_n}(x)=\lim_{y\searrow 0}F(x,y)$$
for $x\in U$ we are done by Theorem 3.1.
\hfill$\blacksquare$

\subsection{Strong quasianalyticity - Proof of Theorem B}

Let $U\subset \IR^n$ be open and let $g:U\to \IR$ be a function. Let $a\in U$. 
The function $g$ is {\bf $N$-flat} at $a$ if $g$ is $C^N$ at $a$ and all partial derivatives of $g$ of order at most $N$
vanish in $a$.
The function $g$ is {\bf flat} at $a$ if $g$ is $C^\infty$ at $a$ and all partial derivatives of $g$ vanish in $a$.
 
\vs{0.5cm}
The asymptotic behaviour of log-analytic function on simple sets implied by the elementary preparation gives the following (see [9] for the corresponding result in polynomially bounded o-minimal structures).

\vs{0.5cm}
{\bf 3.2 Proposition}

\vs{0.1cm}
{\it Let $f:\IR^n\times\IR\to \IR, (x,y)\mapsto f(x,y),$ be a log-analytic function. Then there is $N\in\IN$ such that the following holds for every $x\in \IR^n$:
	If $f(x,-)$ is $N$-flat at $y=0$ then $f(x,-)$ vanishes identically on some open interval around $0\in \IR$.}

\vs{0.1cm}
{\bf Proof:}

\vs{0.1cm}
By also considering the function $f(x,-y)$ it is enough to show that the following holds.
There is $N\in \IN$ such that for every $x\in \IR^n$ with $f(x,-)$ being $N$-flat at $y=0$ we have $f(x,y)=0$ for all $y\in\, ]0,\varepsilon_x[$ for some $\varepsilon_x\in \IR_{>0}$. 

Let $r$ be the log-analytic order of $f$.
By Theorem 2.30 we find a definable cell decomposition $\ma{C}$ of $\IR^n\times\IR$ such that for every simple $C\in\ma{C}$ the cell $C$ is $r$-simple and 
$f|_C$ is purely $r$-log-analytically prepared with respect to $y$.
Let $C\in \ma{C}$ be simple 
and let 
$$\ma{J}=\big(r,\ma{Y}_{r,C}^\mathrm{el},a,q,s,v,b,P\big)$$
be a purely $r$-preparing tuple for $f|_C$.
We show that $a$ vanishes identically.  
Choose $N_C\in \IN$ such that $N_C\geq q_0$.
Let $x\in \pi(C)$.
If
$f(x,-)$ is $N_C$-flat at $y=0$ then
$f(x,-)=o(y^{N_C})$ at $y=0$.
But $|\ma{Y}_{r,C}^\mathrm{el}|^{\otimes q}/y^{N_C}\neq o(1)$ by Remark 2.22.
Hence necessarily $a(x)=0$.

By Remark 2.16 we are done by taking 
$N:=\max\{N_C\mid C\in \ma{C}\mbox{ simple}\}$.
\hfill$\blacksquare$

\vs{0.5cm}
With Proposition 3.2 we can prove Theorem B using familiar connectivity arguments.

\vs{0.5cm}
{\bf Theorem B}

\vs{0.1cm}
{\it Let $U\subset \IR^n$ be a definable domain and let $f:U\to \IR$ be a log-analytic function. Then there is $N\in \IN$ with the following property. If $f$ is $C^N$ and if there is $a\in U$ such that $f$ is $N$-flat at $a$ then $f$ vanishes identically.}

\vs{0.1cm}
{\bf Proof:}

\vs{0.1cm}
Let
$$V:=\big\{(x,z,t,y)\in U\times S^{n-1}\times \IR\times \IR\;\big\vert\; x+(y-t)z\in U\big\}$$
and
$$F:V\times \IR, (x,z,t,y)\mapsto f\big(x+(y-t)z\big).$$
Then $F$ is log-analytic. By Proposition 3.2 there is $N\in \IN$ such that the following holds for every $(x,z,t)\in U\times S^{n-1}\times\IR$ such that $(x,z,t,0)\in V$. If $F(x,z,t,-)$ is $N$-flat at $y=0$ then $F(x,z,t,-)$ vanishes identically on an open interval around $0\in \IR$.

We assume that $f$ is $C^N$. Note that then $F$ is $C^N$.
Let $a\in U$ such that $f$ is $N$-flat at $a$. We show that this implies that $f$ vanishes identically and are done.
We start with the following.

\vs{0.2cm}
{\bf Claim:}
There is $r\in \IR_{>0}$ such that $f$ vanishes identically on $B(a,r)$.

\vs{0.1cm}
{\bf Proof of the claim:}
Let $r\in \IR_{>0}$ be such that $B(a,4r)\subset U$.
Then 
$$W:=B(a,r)\times S^{n-1}\times \,]-r,r[\,\times\, ]-r,r[\,\subset V.$$
Given $z\in S^{n-1}$ we have that $F(a,z,0,-)$ is $N$-flat at $y=0$. 
Then by above $F(a,z,0,-)$ vanishes on some open interval around $0$.
Fix $z\in S^{n-1}$ and let
$A_z$ be the set of all $t\in\, ]-r,r[$ such that $F(a,z,t,-)$ vanishes identically on some open interval around $0$.
Then $A_z\neq \emptyset$ since $0\in A_z$ by above. Clearly $A_z$ is open.
Let $t\in \overline{A_z}\,\cap\, ]-r,r[$. Then $F(a,z,t,-)$ is $N$-flat at $0$. Hence by above $F(a,z,t,-)$ vanishes identically on some open interval around $0$.
Therefore $A_z$ is closed in $]-r,r[$. Since intervals are connected we obtain that $A_z=\,]-r,r[$ and hence that $F(a,z,-,-)$ vanishes identically on $]-r,r[\,\times\,]-r,r[$.
Since $z\in S^{n-1}$ is arbitrary we get that $f$ vanishes identically on $B(a,r)$.
\hfill$\blacksquare_\mathrm{Claim}$

\vs{0.2cm}
Let $X$ be the set of all $x\in U$ such that $f$ vanishes identically on some open ball around $x$.
Then $X\neq \emptyset$ since $a\in X$ by the above. Clearly $X$ is open.
Let $x\in \overline{X}\cap U$. Then $f$ is $N$-flat at $x$. Again by the claim we get that $x\in X$.
Therefore $X$ is closed in $U$. Since $U$ is connected we obtain that $X=U$ and hence that $f$ vanishes identically on $U$.
\hfill$\blacksquare$

\subsection{Tamm's theorem - Proof of Theorem C}

With our result on pure preparation of log-analytic functions on simple sets we can establish the parametric version of Tamm's theorem for log-analytic functions. For this we adapt the reasoning of Van den Dries and Miller [4]. 
The theorem below deals with the case of a parameterized family of unary functions where the adaptions are most extensive.

\vs{0.5cm}
{\bf 3.3 Theorem}

\vs{0.1cm}
{\it Let $f:\IR^n\times\IR\to \IR, (x,y)\mapsto f(x,y),$ be log-analytic. Then there is $N\in \IN$ such that the following holds for all $x\in \IR^n$:
	If $f(x,-)$ is $C^N$ at $0$ then $f(x,-)$ is real analytic at $0$.}

\vs{0.1cm}
{\bf Proof:}

\vs{0.1cm}
Let $r\in \IN_0$ be the log-analytic order of $f$.
By Theorem 2.30 we find a definable cell decomposition $\ma{C}$ of $\IR^n\times\IR$ such that for every simple $C\in\ma{C}$ the cell $C$ is $r$-simple and 
$f|_C$ is purely $r$-log-analytically prepared with respect to $y$.
Fix a simple cell $C\in \ma{C}$ and set $\eta_x:=\sup C_x$ for $x\in \pi(C)$.
Let $\ma{J}=\big(r,\ma{Y}_{r,C}^\mathrm{el},a,q,s,v,b,P\big)$ be a purely $r$-preparing tuple for $g:=f|_C$. 
We set $\ma{Y}:=\ma{Y}_{r,C}^\mathrm{el}$.
Let $\sum_{\alpha\in \IN_0^s}c_\alpha X^\alpha$ be the power series expansion of $v$.  
Let
$$\Gamma_1:=\big\{\alpha\in \IN_0^s\;\big\vert\;^tP\alpha+q\in \IN_0\times \{0\}^r\big\}$$
and
$\Gamma_2:=\IN_0^s\setminus \Gamma_1.$
Set
$v_1:=\sum_{\alpha\in \Gamma_1}c_\alpha X^\alpha$ and $v_2:=\sum_{\alpha\in \Gamma_2}c_\alpha X^\alpha.$
For $l\in \{1,2\}$ let
$$g_l:C\to \IR, (x,y)\mapsto a(x)|\ma{Y}(y)|^{\otimes q}v_l\big(b_1(x)|\ma{Y}(y)|^{\otimes p_1},\ldots,b_s(x)|\ma{Y}(y)|^{\otimes p_s}\big).$$
Then $g_1,g_2$ are log-analytic and $g=g_1+g_2$.
For $k\in \IN_0$
let
$$\Gamma_{1,k}:=\big\{\alpha\in \IN_0^s\;\big\vert\; ^tP\alpha+q=(k,0,\ldots,0)\big\}\subset \Gamma_1$$
and
$$d_k:\pi(C)\to \IR,x\mapsto
a(x)\sum_{\alpha\in \Gamma_k}c_\alpha\prod_{i=1}^sb_i(x)^{\alpha_i}.$$
Then
$$g_1(x,y)=\sum_{k=0}^\infty d_k(x)y^k$$
for $(x,y)\in C$.
The series to the right converges absolutely on $C$ and therefore $g_1$ extends to a well defined extension
$$
\hat{g}_1:\hat{C}:=\big\{(x,y)\in \IR^{n+1}\;\big\vert\; x\in \pi(C), -\eta_x<y<\eta_x\big\}\to \IR,$$
$$(x,y)\mapsto \sum_{k=0}^\infty d_k(x)y^k,$$
such that $\hat{g}_1(x,-)$ is real analytic at $0$ for every $x\in \pi(C)$.

\vs{0.2cm}
{\bf Claim 1:} The function $\hat{g}_1$ is log-analytic. 

\vs{0.1cm}
{\bf Proof of Claim 1:}
Clearly $\hat{C}$ is definable. We show that $\hat{g}_1$ is log-analytic on $\hat{C}\cap (\IR^n\times \IR_{<0})$ and on $\hat{C}\cap (\IR^n\times \{0\})$ and are done. 
For the former
let
$\Gamma_{1,e}:=\bigcup_{k\;\mathrm{even}}\Gamma_{1,k}$ and $\Gamma_{1,o}:=\bigcup_{k\;\mathrm{odd}}\Gamma_{1,k}$. Set $v_{1,e}:=\sum_{\alpha\in \Gamma_{1,e}}c_\alpha X^\alpha$ and $v_{1,0}:=\sum_{\alpha\in \Gamma_{1,o}}c_\alpha X^\alpha$.
Then for $(x,y)\in \hat{C}$ with $y< 0$ we have
$$\hat{g}_1(x,y)=a(x)|\ma{Y}(-y)|^{\otimes q}\Big(v_{1,e}\big(b_1(x)|\ma{Y}(-y)|^{\otimes p_1},\ldots,b_s(x)|\ma{Y}(-y)|^{\otimes p_s}\big)$$
$$-v_{1,o}\big(b_1(x)|\ma{Y}(-y)|^{\otimes p_1},\ldots,b_s(x)|\ma{Y}(-y)|^{\otimes p_s}\big)\Big)$$
which implies the desired assertion.
For the latter note that
$$\hat{g}_1(x,0)=\lim_{y\searrow 0}g_1(x,y)$$ 
for $x\in \pi(C)$. We are done by Theorem 3.1.
\hfill$\blacksquare_{\mathrm{Claim}\,1}$

\vs{0.5cm}
Set $$\hat{g}_2:\hat{C}\to \IR, (x,y)\mapsto f(x,y)-\hat{g}_1(x,y).$$
Then $\hat{g}_2$ is log-analytic by Claim 1 and $\hat{g}_2|_C=g_2$. 
Fix $x^*\in \pi(C)$. 
Let
$$\Lambda:=\big\{^tP\alpha+q\;\big\vert\;\alpha\in \Gamma_2\big\}.$$ 
Then $\Lambda\subset \IQ^{r+1}\setminus (\IN_0\times \{0\}^r)$.
For $\lambda\in \Lambda$, let
$$\Gamma_{2,\lambda}:=\big\{\alpha\in \IN_0^s\;\big\vert\;^tP\alpha+q=\lambda\big\}$$
and
$$e_{x^*,\lambda}:= a(x^*)\sum_{\alpha\in \Gamma_{2,\lambda}}c_\alpha\prod_{i=0}^sb_i(x^*)^{\alpha_i}\in \IR.$$
Then
$$\hat{g}_2(x^*,-)=\sum_{\lambda\in \Lambda}e_{x^*,\lambda}|\ma{Y}|^{\otimes \lambda}$$
on $]0,\eta_{x^*}[$.
Let
$$\Lambda_{x^*}:=\{\lambda\in \Lambda \mid e_{x^*,\lambda}\neq 0\big\}.$$
If $\Lambda_{x^*}= \emptyset$ then $\hat{g}_2(x^*,-)=0$ on $]0,\eta_{x^*}[$. If $\Lambda_{x^*}\neq \emptyset$ there is $\mu_{x^*}=(\mu_{x^*,0},\ldots,\mu_{x^*,r})\in \Lambda_{x^*}$ such that
$|\ma{Y}|^{\otimes \lambda}=o(|\ma{Y}|^{\otimes \mu_{x^*}})$ for all $\lambda\in \Lambda_{x*}$ with $\lambda\neq \mu_{x^*}$.

\vs{0.2cm}
{\bf Claim 2:}
Assume that $\Lambda_{x^*}\neq\emptyset$. Let $M\in \IN$ be such that $f(x^*,-)$ is $C^M$ at $0$.
Then  $\mu_{x^*,0}\geq M$.

\vs{0.1cm}
{\bf Proof of Claim 2:}
Assume that $\mu_{x^*,0}<M$. 

\vs{0.2cm}
{\bf Case 1:} $\mu_{x^*,0}\in \IN_0$.\\
Then $m:=\mu_{x^*,0}+1\leq M$.
Differentiating $g_2$ $m$-times with respect to
$y$ we see with Proposition 2.24 (note that $\mu_{x^*}\neq 0$) that there is $\beta=(-1,\beta_1,\ldots,\beta_r)\in \IQ^{r+1}$ such that
$$\lim_{y \searrow 0}\frac{\partial^mg_2/\partial y^m(x^*,y)}{|\ma{Y}(y)|^{\otimes \beta}}\in \IR^*.$$
Since $\hat{g}_2(x^*,-)=f(x^*,-)-\hat{g}_1(x^*,-)$ is $C^M$ at $0$ we obtain that
$$\lim_{y \searrow 0}\frac{\partial^m\hat{g}_2}{\partial y^m}(x^*,y)=\frac{\partial^m\hat{g}_2}{\partial y^m}(x^*,0)\in \IR,$$
which contradicts that $\hat{g}_2(x^*,-)$ extends $g_2(x^*,-)$.

\vs{0.2cm}
{\bf Case 2:} $\mu_{x^*,0}\notin \IN_0$.\\
Then $m:=\lceil \mu_{x^*,0}\rceil \leq M$. Differentiating $g_2$ $m$-times with respect to
$y$ we see with Proposition 2.24 (note that $\mu_{x^*}\neq 0$) that there is $\beta=(\beta_0,\beta_1,\ldots,\beta_r)\in \IQ^{r+1}$ with $\beta_0<0$ such that
$$\lim_{y \searrow 0}\frac{\partial^mg_2/\partial y^m(x^*,y)}{|\ma{Y}|^{\otimes \beta}}\in \IR^*.$$
But $\hat{g}_2(x^*,-)=f(x^*,-)-\hat{g}_1(x^*,-)$ is $C^M$ at $0$. We get the same contradiction as in Case 1.
\hfill$\blacksquare_{\mathrm{Claim}\,2}$

\vs{0.2cm}
{\bf Claim 3:}
Let $M\in \IN$ be such that $f(x^*,-)$ is $C^M$ at $0$. Then
$\hat{g}_2(x^*,-)$ is $(M-1)$-flat at $0$.

\vs{0.1cm}
{\bf Proof of Claim 3:}

\vs{0.2cm}
{\bf Case 1:} $\Lambda_{x^*}=\emptyset$.\\
Then $\hat{g}_2(x^*,-)=0$ on $]0,\eta_{x^*}[$ and we are clearly done.

\vs{0.2cm}
{\bf Case 2:} $\Lambda_{x^*}\neq \emptyset$.\\
By Claim 1 we obtain that $\mu_{x^*,0}\geq M$.
Hence we obtain by Proposition 2.24 and Remark 2.23 that
$$\lim_{y \searrow 0}\frac{\partial^m g_2}{\partial y^m}(x^*,y)=0$$
for all $m\in \{0,\ldots,M-1\}$.
Since $\hat{g}_2(x^*,-)=f(x^*,-)-\hat{g}_1(x^*,-)$ is $C^M$ at $0$ and since $\hat{g}_2(x^*,-)$ extends $g_2(x^*,-)$ we are done.
\hfill$\blacksquare_{\mathrm{Claim}\,3}$

\vs{0.2cm}
Since $\hat{g}_2$ is log-analytic we find by Proposition 3.2 some $K_C\in \IN$ such that the following holds for every $x\in \pi(C)$. If $\hat{g}_2(x,-)$ is $K_C$-flat at $y=0$ then 
$\hat{g}_2(x,-)$ vanishes identically on some open interval around $0$. 
Set $N_C:=K_C+1$. Assume that $f(x,-)$ is $C^{N_C}$ at $0$. 
Then by Claim 3 $\hat{g}_2(x,-)$ is $K_C$-flat and hence by the above that $f(x,-)=\hat{g}_1(x,-)$ on some open interval around $0$. Since $\hat{g}_1(x,-)$ is real analytic at $0$ we get that $f(x,-)$ is real analytic at $0$. 
By Remark 2.16 we are done by taking
$$N:=\max\{N_C\mid C\in \ma{C}\mbox{ simple}\}.$$
\hfill$\blacksquare$

\vs{0.5cm}
{\bf 3.4 Corollary}

\vs{0.1cm}
{\it Let $f:\IR^n\times\IR\to \IR, (x,y)\mapsto f(x,y),$ be log-analytic such that $f(x,-)$ is real analytic at $0$ for every $x\in \IR^n$. Then there is a definable cell decomposition $\ma{B}$ of $\IR^n$
	such that $B\to \IR, x\mapsto d^k/dy^k f(x,0),$ is real analytic for every $B\in \ma{B}$ and every $k\in \IN_0$.}

\vs{0.1cm}
{\bf Proof:}

\vs{0.1cm}
Using the notation of the previous proof we have $f(x,y)=\hat{g}_1(x,y)$ for all $(x,y)\in C$ where $C$ is a simple cell of the constructed cell decomposition $\ma{C}$. 
Since functions definable in $\IR_{\an,\exp}$ are  piecewise real analytic (see [5, Section 4]) we find a cell decomposition $\ma{D}$ of $\pi(C)\subset \IR^n$ such that the coefficient $a$ and the base functions $b_1,\ldots,b_s$ are real analytic on every $D\in \ma{D}$. 
Hence on each $D\in \ma{D}$ the coefficients $d_k$ of $\hat{g}_1$ are real analytic for every $k\in \IN_0$.
Since $d^k/dy^k f(x,0)=k! d_k(x)$ we are done by Remark 2.16.
\hfill$\blacksquare$

\vs{0.5cm}
Now we can finish as in [4, Section 5].
For the reader's convenience we give the details.

\vs{0.5cm}
Let $U\subset \IR^n$ be open and let $g:U\to \IR$ be a function. Let $a\in U$.
Let $k\in \IN$. Then $g$ is called {\bf $k$-times Gateaux-differentiable} or $G^k$ at $a$ if $y\mapsto g(a+yz)$ is $C^k$ at $y=0$ for every $z\in \IR^n$ and $z\mapsto \big(d^kg(a+yz)/dy^k\big)(0)$ is given by a homogeneous polynomial in $z$ of degree $k$.
The function $g$ is called $G^\infty$ at $a$ if $g$ is $G^k$ at $a$ for every $k\in \IN$. 
The following holds:

\vs{0.5cm}
{\bf 3.5 Fact} ([4, Proposition 2.2])

\vs{0.1cm}
{\it Let $U\subset \IR^n$ be open and let $g:U\to \IR$ be a function.
Let $a\in U$. The following are equivalent.
\begin{itemize}
\item[(i)] The function $g$ is real analytic at $a$.
\item[(ii)] The function $g$ is $G^\infty$ at $a$ and there is $\varepsilon\in \IR_{>0}$ such that for every $z\in \IR^n$ with $|z|\leq 1$ the function $y\mapsto g(a+yz)$ is defined and real  analytic on $]-\varepsilon,\varepsilon[$. 
\end{itemize}}

\vs{0.2cm}
{\bf 3.6 Proposition}

\vs{0.1cm}
{\it Let $f:\IR^n\times\IR^m\to \IR, (x,u)\mapsto f(x,u), $ be a log-analytic function and let $k\in \IN$. Then there is a log-analytic function $w_k:\IR^n\times\IR^m\times\IR^m\to \IR, (x,u,v)\mapsto w_k(x,u,v),$ such that the following is equivalent for every $(x,u)\in \IR^n\times\IR^m$.
\begin{itemize}
	\item[(i)] The function $f(x,-)$ is $G^k$ at $u$.
	\item[(ii)] It is $w_k(x,u,v)=0$ for every $v\in \IR^m$. 
\end{itemize}}
{\bf Proof:}

\vs{0.1cm}
For $k\in \IN$ let
$W_k$ be the set of all $(x,u)\in \IR^n\times\IR^m$ such that 
$y\mapsto f(x,u+yv)$ is $k$-times differentiable in $0$ for every $v\in \IR^m$.
We define
$$\Phi_k:\IR^n\times\IR^m\times\IR^m\to \IR,
(x,u,v)\mapsto \left\{\begin{array}{ccc}
\frac{d^kf(x,u+yv)}{dy^k}(0),&&(x,u)\in W_k,\\
&\mbox{if}&\\
1,&&(x,u)\notin W_k.\\
\end{array}\right.$$
Then $\Phi_k$ is log-analytic by Theorem A.
Let $\nu(k)$ be the dimension of the real vector space of homogeneous real polynomials of degree $k$ in the variables $V:=(V_1,\ldots,V_m)$
and let $M_1(V),\ldots,M_{\nu(k)}(V)$ be the homogeneous monomials of degree $k$ in $V$.
There are points $p_{k,1},\ldots,p_{k,\nu(k)}\in \IR^m$ and linear functions
$a_1,\ldots,a_{\nu(k)}:\IR^{\nu(k)}\to \IR$ such that for all $s:=(s_1,\ldots,s_{\nu(k)})\in \IR^{\nu(k)}$ 
$$P_k(s,V):=\sum_{j=1}^{\nu(k)}a_j(s)M_j(V)\in \IR[V]$$
is the unique homogeneous polynomial of degree $k$ with
$P_k(s,p_{k,i})=s_i$ for all $i\in \{1,\ldots,\nu(k)\}$ (see [4, Section 2.3]).
Set
$$\hat{w}_k:\IR^n\times\IR^m\times\IR^m\to\IR, (x,u,v)\mapsto P_k\big(\Phi_k(x,u,p_{k,1}),\ldots,\Phi_k(x,u,p_{k,\nu(k)}),v\big),$$
and
$w_k:=\hat{w}_k-\Phi_k$. Then $w_k$ is log-analytic. We show that it fulfils the requirements. 

\vs{0.2cm}
$(i)\Rightarrow (ii)$:
Let $(x,u)\in \IR^n\times\IR^m$ be such that $f(x,-)$ is $G^k$ at $u$.
Then $(x,u)\in W_k$ and
$v\mapsto \Phi_k(x,u,v)$ is a homogeneous polynomial of degree $k$.
By the definition of $\hat{w}_k$ we have
$\hat{w}_k(x,u,p_{k,i})=\Phi_k(x,u,p_{k,i})$ for all $i\in \{1,\ldots,\nu(k)\}$.
By the uniqueness of $P_k$ we obtain that
$$\hat{w}_k(x,u,v)=P_k\big(\Phi_k(x,u,p_{k,1}),\ldots,\Phi_k(x,u,p_{k,\nu(k)}),v\big)=\Phi_k(x,u,v)$$
and therefore $w_k(x,u,v)=0$
for all $v\in \IR^m$.

\vs{0.2cm}
$(ii)\Rightarrow (i)$: 
Let $(x,u)\in \IR^n\times\IR^m$ be such that $w_k(x,u,v)=0$ for all $v\in \IR^m$.
Then
$\Phi_k(x,u,v)=\hat{w}_k(x,u,v)$ for all $v\in \IR^m$ and therefore
$v\mapsto \Phi_k(x,u,v)$ is a homogeneous polynomial of degree $k$.
Since $k\geq 1$ it is not constant. 
Hence we get that $(x,u)\in W_k$ and consequently $f$ is $G^k$ at $u$.
\hfill$\blacksquare$

\vs{0.5cm}
{\bf 3.7 Proposition}

\vs{0.1cm}
{\it Let $f:\IR^n\times\IR^m\to \IR, (x,u)\mapsto f(x,u),$ be a log-analytic function.
Then there is $N\in \IN$ such that the following holds for every $(x,u)\in \IR^n\times\IR^m$.
If $f(x,-)$ is $G^N$ at $u$ then $f(x,-)$ is $G^\infty$ at $u$.}

\vs{0.1cm}
{\bf Proof:}

\vs{0.1cm}
By Theorem 3.3 there is $K\in \IN$ such that the following holds for every $(x,u,v)\in \IR^n\times\IR^m\times\IR^m$.
If $y\mapsto f(x,u+yv)$ is $C^K$ at $0$ then $y\mapsto f(x,u+yv)$ is real analytic at $0$.
Letting as in the previous proof $W_k$ for $k\in \IN$ to be the set of all $(x,u)\in \IR^n\times\IR^m$ such that $y\mapsto f(x,u+yv)$ is $k$-times differentiable at $0$ for every $v\in \IR^m$ we get that
$W_k=W_K$ for all $k\geq K$. 
Letting also $\Phi_k$ for $k\in \IN$ be as in the previous proof and dealing with $\Phi_1,\ldots,\Phi_{K-1}$ we find by Corollary 3.4 a definable cell decomposition $\ma{C}$ of $\IR^n\times\IR^m\times\IR^m$ such that $\Phi_k|_C$ is real analytic for every $C\in \ma{C}$ and all
$k\in \IN$. 
Constructing $w_k:\IR^n\times\IR^m\times\IR^m\to \IR, (x,u,v)\to w_k(x,u,v),$ for $k\in \IN$ as in the previous proof we see that $w_k|_C$ is real analytic for every $k\in \IN$. 
By Tougeron [13] (see also [4, Proposition 1.6]) we find for every $C\in\ma{C}$ some $N_C\in \IN$ such that 
$$\bigcap_{k\in \IN}\big\{(x,u,v)\in C\;\big\vert\; w_k(x,u,v)=0\big\}=\bigcap_{k\leq N_C}\big\{(x,u,v)\in C\;\big\vert\; w_k(x,u,v)=0\big\}.$$
Let $N:=\max\{N_C\mid C\in\ma{C}\}$. Then 
$$\bigcap_{k\in \IN}\big\{(x,u,v)\in \IR^n\times\IR^m\times\IR^m\;\big\vert\; w_k(x,u,v)=0\big\}=$$
$$=\bigcap_{k\leq N}\big\{(x,u,v)\in \IR^n\times\IR^m\times\IR^m\;\big\vert\; w_k(x,u,v)=0\big\}.$$
	Hence for every $(x,u)\in \IR^n\times\IR^m$ we have that $f(x,-)$ is $G^\infty$ at $u$ if and only if $w_k(x,u,v)=0$ for all $k\in \IN$ and all $v\in \IR^m$ if and only if $w_k(x,u,v)=0$ for all $k\in \{1,\ldots,N\}$ and all $v\in \IR^m$
	if and only if $f(x,-)$ is $G^N$ at $u$.
	\hfill$\blacksquare$
	
\vs{0.5cm}
{\bf Theorem C} 

\vs{0.1cm}
{\it Let $f:\IR^n\times\IR^m\to \IR, (x,y)\mapsto f(x,y),$ be a log-analytic function.
	Then there is $N\in \IN$ such that the following holds for every $(x,y)\in \IR^n\times\IR^m$.
	If $f(x,-)$ is $C^N$ at $y$ then $f(x,-)$ is real analytic at $y$.} 

\vs{0.1cm}
{\bf Proof:}

\vs{0.1cm}
By Proposition 3.7 there is $N_1\in \IN$ such that if $f(x,-)$ is $G^{N_1}$ at $y$ then $f(x,-)$ is $G^\infty$ at $y$.
Let $F:\IR^n\times\IR^m\times \IR^m\times\IR\to \IR, (x,y,z,t)\mapsto f(x,y+tz).$
By Theorem 3.3 there is $N_2\in \IN$ such that if $F(x,y,z,-)$ is $C^{N_2}$ at $0$ then $F(x,y,z,-)$ is real analytic at $0$.
Taking $N:=\max\{N_1,N_2\}$ we are done by Fact 3.5.
\hfill$\blacksquare$

\vs{0.5cm}
{\bf 3.8 Corollary}

\vs{0.1cm}
{\it Let $f:\IR^n\times\IR^m\to \IR, (x,y)\mapsto f(x,y),$ be a log-analytic function.
	Then the set of all $(x,y)\in \IR^n\times\IR^m$ such that $f(x,-)$ is real analytic at $y$ is definable.}

\vs{1cm}
\newpage
\noi \footnotesize{\centerline{\bf References}
	\begin{itemize}
		\item[(1)] 
		R. Cluckers and D. Miller:
		Stability under integration of sums of products of real globally subanalytic functions and their logarithms.
		{\it Duke Math. J.} {\bf 156} (2011), no. 2, 311-348.
		\item[(2)]
		L. van den Dries: Tame Topology and O-minimal Structures. {\it London Math. Soc. Lecture Notes Series} {\bf 248}, Cambridge University Press, 1998.
		\item[(3)] 
		L. van den Dries, A. Macintyre, D. Marker:
		The elementary theory of restricted analytic fields with exponentiation.
		{\it Annals of Mathematics} {\bf 140} (1994), 183-205.
	    \item[(4)]
		L. van den Dries, C. Miller:
		Extending Tamm's theorem. {\it Ann. Inst. Fourier} {\bf 44}, no. 5 (1994), 1367-1395.
		\item[(5)]
		L. van den Dries and C. Miller:
		Geometric categories and o-minimal structures.
		{\it Duke Math. J.} {\bf 84} (1996), no. 2, 497-540.
		\item[(6)]
		L. van den Dries and P. Speissegger:
		O-minimal preparation theorems.
		Model theory and applications, 87-116, Quad. Mat., 11, Aracne, Rome, 2002.
		\item[(7)]
		T. Kaiser: Integrations of semialgebraic functions and integrated Nash functions. {\it Mathematische Zeitschrift} {\bf 275} (2013),
		349-366. 
		\item[(8)]
		J.-M. Lion, J.-P. Rolin:
		Th\'{e}or\`{e}me de pr\'{e}paration pour les fonctions logarithmico-exponentielles.
		{\it Ann. Inst. Fourier} {\bf 47}, no. 3 (1997), 859-884.
		\item[(9)]
		J.-M. Lion, J.-P. Rolin:
		Int\'egration des fonctions sous-analytices et volumes des sous-ensembles sous-analytiques. 
		{\it Ann. Inst. Fourier} {\bf 48}, no. 3 (1998), 755-767.
		\item[(10)] 
		C. Miller: 
		Infinite differentiability in polynomially bounded o-minimal structures. {\it Proc. Amer. Math. Soc.} {\bf 123} (1995), 2552-2555.
		\item[(11)]
		W. Paw\l ucki, A. Pi\c{e}kosz: A remark on the Lion-Rolin Preparation Theorem for LA-functions. {\it Ann. Pol. Math.}, no. 2 (1999), 195-197.
        \item[(12)]
        M. Tamm: Subanalytic sets in the calculus of variations. {\it Acta Mathematica} {\bf 146} (1981), 167-199.
        \item[(13)] 
		J.-C. Tougeron: Alg\`{e}bres analytiques topologiquement noeth\'{e}riennes. Th\'{e}orie de Khovanskii. {\it Ann. Inst. Fourier (Grenoble)} {\bf 41}, no. 4 (1991), 823-840.
\end{itemize}}

\vs{0.5cm}
Tobias Kaiser\\
University of Passau\\
Faculty of Computer Science and Mathematics\\
tobias.kaiser@uni-passau.de\\
D-94030 Germany

\vs{0.5cm}
Andre Opris\\
University of Passau\\
Faculty of Computer Science and Mathematics\\
andre.opris@uni-passau.de\\
D-94030 Germany
\end{document}